\documentclass[11pt]{article}

\usepackage[sort]{cite}       
\usepackage{amssymb}          
\usepackage{latexsym}         
\usepackage{amsmath}          
\usepackage[latin1]{inputenc} 
\usepackage{eucal}            
\usepackage{bbm}              
\usepackage{theorem}          
\usepackage{stmaryrd}         
\usepackage{mathrsfs}         

\title{A $C^*$-Algebraic Model for Locally Noncommutative Spacetimes}

\author{\textbf{Jakob G. Heller}\thanks{Jakob.Heller@physik.uni-freiburg.de},
  \addtocounter{footnote}{2}
  \textbf{Nikolai Neumaier}\thanks{Nikolai.Neumaier@physik.uni-freiburg.de},
  \addtocounter{footnote}{0}
  \textbf{Stefan Waldmann}\thanks{Stefan.Waldmann@physik.uni-freiburg.de}
  \\[0.1cm]
  Fakult{\"a}t f{\"u}r Mathematik und Physik\\
  Albert-Ludwigs-Universit{\"a}t Freiburg\\
  Physikalisches Institut\\
  Hermann Herder Stra{\ss}e 3\\
  D 79104 Freiburg\\
  Germany
}

\date{September 2006\\[0.5cm] FR-THEP 2006/14}

\renewcommand{\mathbb}[1]{\mathbbm{#1}} 

\newcommand{\id}         {\operatorname{\mathsf{id}}}   
\newcommand{\supp}       {\operatorname{\mathrm{supp}}}

\newcommand{\SP}[1]      {\left\langle{#1}\right\rangle}

\newcommand{\abs}[1]    {\left|{#1}\right|} 
\newcommand{\norm}[1]    {\left\|{#1}\right\|}            
 
\newcommand{\I}          {\mathrm{i}}
\newcommand{\E}          {\mathrm{e}}
\newcommand{\D}          {\operatorname{\mathrm{d}}}

\newcommand{\Ver}        {\operatorname{\mathrm{Ver}}}

\newcommand{\Anti}       {\Lambda}

\newcommand{\ver}        {\mathsf{ver}}
\newcommand{\Schouten}[1]{\left\llbracket{#1}\right\rrbracket}

\newcommand{\starp}      {\mathbin{\star_p}}
\newcommand{\tstar}      {\mathbin{\tilde{\star}}}
\newcommand{\tstarp}     {\mathbin{\tilde{\star}_p}}

\newtheorem{lemma} {Lemma} [section]
\newtheorem{proposition} [lemma] {Proposition}

\newtheorem{definition}[lemma] {Definition}

\theorembodyfont{\rmfamily}

\newtheorem{remark}[lemma]{Remark}

\newenvironment{proof}[1][{}]{
  \par\noindent
  \textsc{Proof{#1}:}
}
{
  \hspace*{\fill} $\blacksquare$\newline
}

\numberwithin{equation}{section}

\begin{document}

\maketitle

\begin{abstract}
    Locally noncommutative spacetimes provide a refined notion of
    noncommutative spacetimes where the noncommutativity is present
    only for small distances. Here we discuss a non-perturbative
    approach based on Rieffel's strict deformation quantization.  To
    this end, we extend the usual $C^*$-algebraic results to a
    pro-$C^*$-algebraic framework.
\end{abstract}

%
%

\section{Introduction}
\label{sec:Intro}

In the last years, models for spacetime at small distances using
noncommutative geometric structures in the sense of Connes
\cite{connes:1994a} have received increasing attention, though the
ideas of making spacetime noncommutative can be traced back in history
quite far, see e.g.~\cite{rieffel:1997a,
  doplicher.fredenhagen.roberts:1995a, jurco.et.al.:2001a,
  jurco.schupp.wess:2000a, connes.douglas.schwarz:1998a,
  bahns.doplicher.fredenhagen.piacitelli:2003a} for some recent works.
However, most of these models have severe difficulties, both from the
technical but also from the conceptual point of view: in particular,
the noncommutative Minkowski spacetime with a constant Poisson
structure $\theta$ and hence the usual Weyl-Moyal star product is
certainly not `geometric' at all but refers to a very particular
symmetry, namely of a flat spacetime. Moreover, the fact that $\theta$
is constant leads to global effects which should allow to observe the
noncommutativity already at macroscopic distances. The famous UV/IR
mixing can be seen as an indication for this. In some sense, this has
to be expected to be a generic effect as long as the support of
$\theta$ coincides with the whole spacetime.

In \cite{bahns.waldmann:2006a} a slightly different approach to
\emph{locally} noncommutative spacetimes was established by taking
seriously the wish that noncommutativity should only be visible
\emph{at small distances}. This naturally leads to a framework where
not the spacetime $M$ but $M \times M$ is equipped with a star product
$\tstar$ as one needs two points in order to speak of `distance'.
Then, $\tstar$ is chosen such that it is only non-trivial close to the
\emph{diagonal} $\Delta_M \subseteq M \times M$. We shall briefly
recall the basic constructions and refer to
\cite{bahns.waldmann:2006a} for the physical motivation and the
interpretation of the model.

Let $(M, \nabla)$ be a manifold with a torsion-free connection, e.g.\ a
(pseudo-)Riemannian manifold $(M, g)$ with the Levi-Civita
connection. We choose an open neighbourhood
$\mathcal{U} \subseteq TM$ of the zero section together with an open
neighbourhood $\mathcal{V} \subseteq M \times M$ of the diagonal such
that the map
\begin{equation}
    \label{eq:PhiDiffeo}
    \Phi: \mathcal{U} \ni v_p 
    \; \mapsto \; 
    \Phi(v_p) = (\exp_p(-v_p), \exp_p(v_p)) \in \mathcal{V}
\end{equation}
is a well-defined diffeomorphism, where $\exp$ is the exponential map
of $\nabla$. We use $\Phi$ to establish geodesic relative coordinates
on $M \times M$ close around $\Delta_M$. In a next step, one requires
a Poisson structure $\theta$ on $TM$ with the following properties:
$\theta \in \Gamma^\infty(\Anti^2 \Ver(TM))$ is \emph{vertical} and
$\supp \theta \subseteq \mathcal{U}$. Thus $\Phi_* \theta =
\tilde{\theta} \in \Gamma^\infty(\Anti^2 T(M \times M))$ is a
well-defined Poisson structure on $M \times M$ which is only
non-trivial in direction of the geodesic relative coordinates.
Finally, one requires $\supp \theta \cap T_pM$ to be compact for all
$p \in M$ which encodes the idea of noncommutativity at small
distances. It was shown in \cite{bahns.waldmann:2006a} that such a
$\theta$ can be quantized by a formal star product $\star$ on $TM$
sharing essentially the same properties: $\star$ is vertical and only
non-trivial in $\mathcal{U}$. Moreover, $\star$ restricts to a star
product $\starp$ on $T_pM$ for each $p \in M$ such that the higher
order cochains of $\starp$ have compact support around $0_p$ in $T_pM
\cap \mathcal{U}$. Hence, $\star$ can be pushed forward to $\tstar$ on
$M \times M$ by $\Phi$. The pull-back $\Phi^*$, defined at least on
functions with support in $\mathcal{V}$, becomes an algebra morphism
between $\tstar$ and $\star$. Thanks to the support properties of
$\starp$ we can push forward $\starp$ by $\exp_p$ to $M$ yielding a
star product $\tstarp$ for each $p \in M$ which is only non-trivial in
a small neighbourhood of $p$. This way, every point obtains its own
small noncommutative neighbourhood.

The aim of this letter is to show that the whole construction can
still be done in a $C^*$-algebraic and also in a pro-$C^*$-algebraic
framework if we restrict ourselves to a particular but still very rich
class of Poisson structures.

This is important and interesting for several reasons: on the one
hand, to set up reasonable quantum field theories on noncommutative
spacetimes a more analytic framework than just formal star products is
needed and $C^*$-algebras have shown to be a good choice
\cite{doplicher.fredenhagen.roberts:1995a,
  bahns.doplicher.fredenhagen.piacitelli:2003a}.  On the other hand,
we believe that our construction, which is a slight variation of
Rieffel's strict deformation quantization by actions of $\mathbb{R}^d$
\cite{rieffel:1993a}, is interesting for its own, independently of the
possible interpretation in the context of noncommutative spacetimes.
In particular, the class of Poisson structures we consider does not
seem to be accessible by other techniques than Rieffel's quantization
like e.g.\ deforming by using generators and relations: in the
geometric framework we intend to work in there are no reasonable
generators and relations neither for the Poisson algebra nor for the
deformation. Finally, the extension to the pro-$C^*$-algebraic world
encodes in an appropriate way the non-compactness of the tangent
bundle. It seems to be an interesting generalization of Rieffel's
original work \cite{rieffel:1993a}, independently of our application
to locally noncommutative spacetimes.

The letter is organized as follows: in
Section~\ref{sec:VerticalPoissonStructures} we show that there is a
large class of vertical Poisson structures meeting all the properties
needed for Rieffel's quantization. Further, we show in the following
section how to apply Rieffel's construction of strict deformation
quantization \cite{rieffel:1993a} to the $C^*$-algebra of bounded
continuous functions $C^0_b(N)$ on a manifold $N$ and to the
pro-$C^*$-algebra of continuous functions $C^0(N)$. The latter for the
particular case of an action which is trivial outside a compact
subset. In Section~\ref{sec:LocallyNoncommutativeSpaceTime} we use the
deformed function algebras to construct a locally noncommutative
spacetime. For the convenience of the reader the appendix gives a
short survey on Rieffel's strict deformation quantization and its
extension to pro-$C^*$-algebras.

%
%

\medskip

\noindent
\textbf{Acknowledgement:} We would like to thank Dorothea Bahns and
Marc Rieffel for valuable discussions and Alan Weinstein for
suggesting the use of the diffeomorphism as in
Lemma~\ref{lemma:NiceDiffeo}.

%
%

\section{Vertical Poisson Structures from Actions of $\mathbb{R}^{2d}$}
\label{sec:VerticalPoissonStructures}

In this section we shall construct vertical Poisson structures on a
real vector bundle $\pi: E \longrightarrow M$ of fibre dimension $n$
arising from a particular action of $\mathbb{R}^{2d}$ for sufficiently
large $d$.  Recall that a $k$-vector field $X \in
\Gamma^\infty(\Anti^k TE)$ is called \emph{vertical} if $X\in
\Gamma^\infty(\Anti^k \Ver(TE))$, where $\Ver(TE) = \ker T\pi
\subseteq TE$ is the vertical subbundle of $TE$. Moreover, if $s \in
\Gamma^\infty(E)$ is a section then the vertical lift $s^\ver \in
\Gamma^\infty(\Ver(TE))$ is defined by $s^\ver(e) = \frac{\D}{\D
  t}\big|_{t=0} (e + t s(\pi(e)))$. We extend the vertical lift to
arbitrary sections of tensor powers of $E$.

We start with the following technical lemma on a particular
diffeomorphism $\Psi: B_1(0) \longrightarrow \mathbb{R}^n$, the use of
which was suggested to us by Alan Weinstein \cite{weinstein:2006a:misc}.
\begin{lemma}
    \label{lemma:NiceDiffeo}
    There exists a diffeomorphism $\Psi: B_1(0) \longrightarrow
    \mathbb{R}^n$ with the following properties:
    \begin{enumerate}
    \item $\Psi$ is the identity on $B_{\frac{1}{2}}(0)$.
    \item $\Psi^* \frac{\partial}{\partial x^i} = X_i$ extends from
        $B_1(0)$ to a smooth vector field on $\mathbb{R}^n$ by setting
        $X_i$ equal to $0$ outside of $B_1(0)$ for $i = 1, \ldots, n$.
    \item $\Psi$ is $\mathrm{O}(n)$-equivariant with respect to the
        canonical action of $\mathrm{O}(n)$ on $B_1(0)$ and
        $\mathbb{R}^n$.
    \end{enumerate}
\end{lemma}
\begin{proof}
    Define $\psi: [0, 1) \longrightarrow [0,
    +\infty)$ by $\psi(t) = t \chi(t) + (1-\chi(t)) \E^{\frac{1}{1 -
      t}}$, where $\chi$ is a suitable cut-off function such that
    $\chi(t) = 1$ for $t \le \frac{1}{2}$ and $\chi(t) = 0$ for $t \ge
    \frac{3}{4}$ and such that $\psi$ is a diffeomorphism. Obviously,
    $\chi$ can be arranged in such a way. Then $\Psi: B_1(0)
    \longrightarrow \mathbb{R}^n$ defined by 
    \[
    \Psi(x) = \frac{x}{\abs{x}} \psi(\abs{x})
    \]
    will fulfill the assertions, where $\abs{x}$ denotes the Euklidian
    norm of $x$.
\end{proof}

Now consider an arbitrary positive definite fibre metric $h$ on $E$
and denote by $B_1^h(0) \subseteq E$ the bundle of open balls of
radius $1$ with respect to $h$.
\begin{lemma}
    \label{lemma:VectorBundleShrinkingDiffeo}
    Let $e_i \in \Gamma^\infty(E)$ with $i \in I$ be a collection of
    sections of $E$. Then there exist vector fields $X_i \in
    \Gamma^\infty(TE)$ on $E$ with the following properties:
    \begin{enumerate}
    \item $X_i \in \Gamma^\infty(\Ver(TE))$ is vertical for all $i \in
        I$.
    \item $[X_i, X_j] = 0$ for all $i, j \in I$.
    \item $X_i$ coincides with the vertical lift $e_i^\ver$ on
        $B_{\frac{1}{2}}^h(0)$ for all $i \in I$.
    \item $\supp X_i \subseteq B_1^h(0)$ for all $i \in I$.
    \end{enumerate}
\end{lemma}
\begin{proof}
    First note that the vertical lifts $e_i^\ver \in
    \Gamma^\infty(\Ver(TE))$ of sections $e_i \in \Gamma^\infty(E)$
    satisfy $[e_i^\ver, e_j^\ver] = 0$ for all $i, j \in I$. Moreover,
    $e_i^\ver$ is constant in fibre directions. As we have chosen a
    metric $h$ we can define the diffeomorphism $\Psi$ from
    Lemma~\ref{lemma:NiceDiffeo} fibrewise with respect to $h$. This
    is well-defined thanks to the $\mathrm{O}(n)$-equivariance. Then
    $X_i = \Psi^* e_i^\ver$ has the desired properties.
\end{proof}

Now we can prove the existence of `many' vertical Poisson structures
with compact support in fibre directions:
\begin{proposition}
    \label{proposition:NicePoissonStructures}
    Let $\gamma \in \Gamma^\infty(\Anti^2E)$ and let $\mathcal{U}
    \subseteq E$ be an open neighbourhood of the zero section of $E$.
    Then there exist $2d$ vector fields $X_1, \ldots, X_d$, $Y_1,
    \ldots, Y_d \in \Gamma^\infty(TE)$ and an open neighbourhood
    $\tilde{\mathcal{U}} \subseteq \mathcal{U}$ of the zero section
    such that:
    \begin{enumerate}
    \item $X_i, Y_i \in \Gamma^\infty(\Ver(TE))$ are vertical for all
        $i = 1, \ldots, d$.
    \item $\supp X_i, \supp Y_i \subseteq \mathcal{U}$ for all $i = 1,
        \ldots, d$ and $\supp X_i \cap E_p$, $\supp Y_i \cap E_p$ are
        compact for all $i = 1, \ldots, d$.
    \item $[X_i, X_j] = [X_i, Y_j] = [Y_i, Y_j] = 0$ for all $i, j =
        1, \ldots, d$.
    \item Each of the vector fields $X_i$, $Y_j$ has complete flow.
    \item $\theta = \sum_{i=1}^d X_i \wedge Y_i$ is a vertical Poisson
        structure such that $\theta$ coincides with the vertical lift
        $\gamma^\ver \in \Gamma^\infty(\Anti^2 \Ver(TE))$ on
        $\tilde{\mathcal{U}}$.
    \end{enumerate}
\end{proposition}
\begin{proof}
    By the Serre-Swan theorem we know that $\Gamma^\infty(E)$ is a
    finitely generated projective module over $C^\infty(M)$, even if
    $M$ is non-compact. Thus there exist sections $e_i \in
    \Gamma^\infty(E)$ and $f^i \in \Gamma^\infty(E^*)$ with $i = 1,
    \ldots, d$, where $d$ is sufficiently large, such that for all $s
    \in \Gamma^\infty(E)$ we have $s = \sum_{i=1}^d e_i f^i(s)$. In
    other words, the sections $e_i$ and $f^i$ form a finite dual
    basis, see e.g. \cite[Lemma 2.9]{lam:1999a} as well as
    \cite[Proposition 4.2]{wells:1980a}. It follows that there exist
    globally defined functions $\gamma^{ij} = - \gamma^{ji} \in
    C^\infty(M)$ such that $\gamma = \frac{1}{2} \sum_{i, j}
    \gamma^{ij} e_i \wedge e_j$. Note that in general $d > n$ (unless
    $E$ is trivial) and hence the $\{e_i\}$ are \emph{not} linearly
    independent. Next, we choose a fibre metric $h$ such that $B_1^h
    (0) \subseteq \mathcal{U}$. Then
    Lemma~\ref{lemma:VectorBundleShrinkingDiffeo} provides us with $d$
    commuting vector fields $X_i$ that coincide with $e_i^\ver$ on
    $B_{\frac{1}{2}}^h(0)$ and satisfy $\supp X_i \subseteq B_1^h(0)$.
    Now define $Y_i = \frac{1}{2} \sum_{j=1}^d \pi^* \gamma^{ij} X_j$
    and set $\theta = \sum_{i=1}^d X_i \wedge Y_i$. Since the $X_i$
    are vertical and the functions $\pi^* \gamma^{ij}$ are constant
    along the fibres it follows that the vector fields $X_i$ and $Y_j$
    still commute and hence $\Schouten{\theta, \theta} = 0$. Thus they
    satisfy the first three requirements. Moreover, since the support
    of the vertical vector fields $X_i$ and $Y_i$ is compact in fibre
    directions, their flows are complete.  Finally, setting
    $\tilde{\mathcal{U}} = B_{\frac{1}{2}}^h (0)$, the Poisson
    structure $\theta$ satisfies the last requirement.
\end{proof}

Conversely, assume we have a smooth action of $\mathbb{R}^{d}$ on $E$
by vertical diffeomorphisms. Then the fundamental vector fields $X_i$,
$i = 1, \ldots, d$, of the action are vertical and $\theta =
\frac{1}{2} \sum_{i=1}^d \Theta^{ij} X_i \wedge X_j$ defines a Poisson
structure on $E$, for all choices of constants $\Theta^{ij} = -
\Theta^{ji} \in \mathbb{R}$. The support condition on $\theta$ simply
means that outside an open neighbourhood of the zero section of $E$
all points are fixed points for the action.
\begin{definition}
    \label{definition:ThetaAdmissible}
    Let $\mathcal{U} \subseteq E$ be an open neighbourhood of the zero
    section and let $\theta \in \Gamma^\infty(\Anti^2 TE)$ be a
    vertical Poisson structure. Then $\theta$ is called
    $\mathcal{U}$-admissible if there exists a smooth action of
    $\mathbb{R}^{d}$ on $E$ by vertical diffeomorphisms and constants
    $\Theta^{ij} = - \Theta^{ji} \in \mathbb{R}$ such that
    \begin{equation}
        \label{eq:ThetaXiXj}
        \theta = \frac{1}{2} \sum_{i, j=1}^d \Theta^{ij} X_i \wedge X_j,
    \end{equation}
    with $\supp X_i \subseteq \mathcal{U}$, where the $X_i$ are the
    (vertical) fundamental vector fields of the action and $\supp X_i
    \cap T_pM$ is compact for all $p \in M$.
\end{definition}
\begin{remark}
    \label{remark:Nontrivial}
    Proposition~\ref{proposition:NicePoissonStructures} says that for
    any $\mathcal{U}$ there exist `many' non-trivial vertical Poisson
    structures which are $\mathcal{U}$-admissible.
\end{remark}

%
%

\section{Strict Deformation Quantization for Actions with Compact
  Support}
\label{sec:RieffelCompactSupport}

In this section, we apply the strict deformation quantization
indroduced in \cite{rieffel:1993a} to construct a noncommutative
product for the $C^*$-algebra of all bounded continuous functions
$C^0_b(N)$ on a manifold $N$ and the pro-$C^*$-algebra of
continuous functions $C^0(N)$ on $N$, respectively.

We choose $d$ vector fields $X_1, \ldots, X_d \in
\Gamma^\infty(TN)$ with the following properties: 
\begin{enumerate}
\item $\supp X_i \subseteq K \subseteq N $ for $i = 1, \ldots, d$,
    where $K$ is a compact subset of $N$.
\item $[X_i, X_j] =0$ for $i,j = 1,\ldots, d$.
\end{enumerate}
The existence of such vector fields which are even non-trivial on $K$
is guaranteed by a slight modification of
Lemma~\ref{lemma:VectorBundleShrinkingDiffeo}. Their flows
$\phi^{X_1}, \ldots, \phi^{X_d}$ determine an action of $\mathbb{R}^d$
on $C^0_b(N)$ by
\begin{equation}
    \alpha(v,f) = \alpha_v (f) = f \circ \phi^{X_1}_{v_1} \circ \cdots
    \circ \phi_{v_d}^{X_d}.
\end{equation}
Due to the properties of $X_1, \ldots, X_d \in \Gamma^\infty(TN)$,
this indeed defines an action. It can easily be seen that the mapping
$v \mapsto \alpha_v (f)$ is continuous for all $f \in C_b^0(N)$.
Indeed we have $\alpha_v (f) = f$ on $N \setminus K$ and as $f$ is
uniformly continuous on $K$, the assertion follows immediately.
Moreover, it is obviously isometric with respect to the supremum norm
on $C_b^0(N)$. Therefore it is possible to apply Rieffel's
construction in this framework:
\begin{definition}
    \label{definition:product}
    Let $\Theta$ be a linear and skew-symmetric operator on
    $\mathbb{R}^d$ with respect to the standard inner product.  The
    noncommutative product $\star: C_b^0(N)^\infty \times
    C_b^0(N)^\infty \longrightarrow C_b^0(N)^\infty$ is defined by
    \begin{equation}
        f \star g 
        = \iint \alpha_{\Theta u}(f) \alpha_v (g) 
        \E^{2 \pi \I u\cdot v} \D u \D v.
    \end{equation}
    Here $C_b^0(N)^\infty$ denotes the space of smooth vectors of
    $\alpha$ in $C_b^0(N)$ (see Appendix~\ref{subsec:deformed}).
\end{definition}
As $C^0_b(N)$ is a $C^*$-Algebra, it is possible to define a norm
$\norm{\,.\,}_\Theta $ on $C_b^0(N)^\infty$ such that the completion
of $(C_b^0(N)^\infty, \star, \norm{\,.\,}_\Theta )$ becomes a
$C^*$-algebra (see Appendix~\ref{subsec:cstar}). Since the flows of
$X_1, \ldots, X_d$ are smooth, it follows that $C_0^\infty(N)
\subseteq C_b^0(N)^\infty$.

Before we will use this construction of a noncommutative product of
$C^0_b(N)$ in the next section to construct locally noncommutative
spacetimes, we present some first results concerning the properties of
the algebra $(C_b^0(N)^\infty, \star, \norm{\,.\,}_\Theta )$.

First we show, that the subalgebra of those functions in
$C^0_b(N)^\infty$ with compact support $C^0_b(N)^\infty \cap C^0_0(N)$
remains a subalgebra with respect to the deformed product.
\begin{proposition}
    \label{proposition:compactsupport}
    The functions in $C^0_b(N)^\infty$ with compact support form a
    subalgebra with respect to the deformed product $\star$, i.e.\ for
    $f,g \in C^0_b(N)^\infty \cap C^0_0(N)$ we have $f \star g \in
    C^0_b(N)^\infty \cap C^0_0(N)$. More explicitly
    \begin{equation}
        \supp (f\star g) \subset (\supp f \cap \supp g) \cup K.
    \end{equation}
\end{proposition}
\begin{proof}
    Let $\chi \in C^\infty_0(N)$ be a cut-off function with $\chi|_K
    \equiv 1$. Then we have $f = \chi f + (1-\chi)f$. Due to $\supp
    (1- \chi)f \cap K = \emptyset$ the function $(1 - \chi)f$ is a
    fixpoint of the action $\alpha$. Thus
    Proposition~\ref{proposition:fixpoint} yields $((1 - \chi)f) \star
    g = (1 -\chi)fg$ and therefore
    \[
        \supp(((1- \chi)f) \star g) \subseteq (N\backslash K) \cap
        \supp f \cap \supp g. 
    \]
    Moreover, we have $\alpha_v (\chi) = \chi$. Therefore, due to
    Proposition~\ref{proposition:fixpoint} we have $\chi f = \chi
    \star f$. The associativity of the deformed product gives $(\chi
    f) \star g = (\chi \star f) \star g = \chi \star (f \star g) =
    \chi (f \star g)$.  Thus, we have the following inclusion for the
    support of $(\chi f) \star g $
    \[
        \supp \left((\chi f) \star  g \right) \subseteq \supp \chi \cap
        \supp(f \star g) \subseteq \supp \chi.
    \]
    Consequently, we have $f \star g = (\chi f + (1- \chi)f) \star g =
    (\chi f) \star g + ((1-\chi)f)\star g$ whence
    \[
        \supp (f\star  g) \subseteq \supp \left( (\chi f) \star  g
        \right) \cup \supp
        \left( ((1-\chi)f)\star  g \right) \subseteq (\supp f \cap
        \supp g) \cup \supp \chi.
    \]
    As for each $q \notin K$, there is a cut-off function $\chi \in
    C^\infty(N)$ with $\chi|_K \equiv 1$ and $\chi(q) = 0$, the
    assertion follows.
\end{proof}

From the proof of the last proposition one immediately sees that the
functions $f \in C^0_b(N)^\infty$ whose supports have an empty
intersection with $K$ form a central $^*$-ideal of $(C_b^0(N)^\infty,
\star, \norm{\,.\,}_\Theta )$.
\begin{proposition}
    \label{proposition:Ideal}
    For $f \in C^0_b(N)^\infty$ with $\supp f \cap K = \emptyset$ we
    have
    \begin{equation}
        f \star g = fg = g \star f 
        \qquad \forall g \in C_b^0(N)^\infty.  
    \end{equation}
    In particular, $\supp(f \star g) \cap K = \emptyset = \supp(g
    \star f) \cap K$.
\end{proposition}
Concerning the states of $(C_b^0(N)^\infty, \star,
\norm{\,.\,}_\Theta)$, we have no general results yet. However, the
following proposition shows that in the situation, where the vector
fields $X_i$ have common zeros, some states of the deformed algebra
coincide with those of the undeformed algebra:
\begin{proposition}
    \label{proposition:DeltaFunktional}
    Let $q \in K$ with $X_i(q) = 0$ for all $i= 1, \ldots, d$. Then we
    have 
    \begin{equation}\label{eq:DeltaFunktional}
        \delta_q(f\star g) = f(q) g(q) = \delta_q(fg) 
        \qquad\forall f,g\in C_b^0(N)^\infty,
    \end{equation}
    where $\delta_q: C^0_b(N) \longrightarrow \mathbb{R}$ denotes the
    $\delta$-functional. For $q\in N \setminus K$ the validity of
    Equation~\eqref{eq:DeltaFunktional} is obvious from the properties
    of the action.
\end{proposition}
\begin{proof}
    For $X_i(q) = 0$ we have $\phi^{X_i}_v(q) = q$. According to the
    presumptions, this holds for all $i = 1, \ldots, d$, such that
    $\alpha_v(f) (q) = f(q)$. Therefore we have with \cite[Corollary
    1.12]{rieffel:1993a}:
    \begin{equation*}
        \delta_q(f \star g)
        = \iint \delta_q(\alpha_{\Theta u}(f)) \delta_q(\alpha_v(g))
        \E^{2 \pi \I u\cdot v} \D u \D v 
        = \iint \delta_q(f)\delta_q(g) \E^{2 \pi \I u \cdot v} \D u \D v 
        = \delta_q(fg). 
    \end{equation*}
\end{proof}

To conclude this section, we want to emphasize that there are at least
two possible variations of the construction above:
\begin{remark}
    \label{remark:CNullImUnendlichen}
    We can use the $C^*$-algebra $C^0_\infty(N)$ of continuous
    functions vanishing at infinity instead of $C^0_b(N)$, i.e.\ the
    $C^*$-completion of $C^0_0(N)$. Then all results remain true for
    the corresponding deformation $(C^0_\infty(N)^\infty, \star,
    \norm{\,.\,}_\Theta)$ of $C^0_\infty(N)$.
\end{remark}
\begin{remark}
    \label{remark:ProCstar}
    Furthermore, we can replace the bounded continuous functions by
    the pro-$C^*$-algebra $C^0(N)$ of \emph{all} continuous functions.
    Here we need the property that the supports of the vector fields
    $X_1, \ldots, X_d$ are compact. Indeed, the action $\alpha$ on
    $C^0(N)$ is strongly continuous since on a compact subset a
    continuous function is uniformly continuous.  Moreover, $\alpha$
    is cofinally isometric, see
    Definition~\ref{definition:CofinaleWirkung}, as for each compactum
    $L \supseteq K$ we clearly have 
    \begin{equation}
        \label{eq:Isometric}
        \norm{f}_L = \norm{\alpha_v(f)}_L \quad \forall v \in
        \mathbb{R}^d
    \end{equation}
    for the sup-norm $\norm{\,.\,}_L$ over $L$.  Hence in the same way
    as in the case of $C^0_b(N)$ we obtain a deformed algebra
    $(C^0(N)^\infty, \star)$, which can be made a
    pre-pro-$C^*$-algebra due to the results of
    Appendix~\ref{subsec:procstar}.  Then a straightforward
    verification shows that the assertions of the
    Propositions~\ref{proposition:compactsupport},
    \ref{proposition:Ideal}, and \ref{proposition:DeltaFunktional}
    literally hold true for $C^0(N)$ in the place of $C^0_b(N)$. Note
    that for \eqref{eq:Isometric} the condition $\supp X_i \subseteq
    K$ is crucial.
\end{remark}

%
%

\section{Locally Noncommutative Spacetimes}
\label{sec:LocallyNoncommutativeSpaceTime}

This section will be devoted to the construction of deformed products
that incorporate the idea of a locally noncommutative spacetime in the
framework of (pro-) $C^*$-algebras, using Rieffel's strict deformation
quantization presented in the preceding section. We will proceed in
several steps first constructing a suitable action of $\mathbb{R}^d$
on $TM$ for an arbitrary smooth manifold $M$ that induces a deformed
product on $C^0(TM)^\infty$ and then using the exponential map of a
connection in $TM$ to obtain induced products on $C^0(M \times
M)^\infty$ and $C^0(M)^\infty$.  Furthermore, we will clarify the
relations between the different products. Clearly, all constructions
will have their equivalents for the cases $C^0_\infty$ and $C^0_b$.

Let $\mathcal{U}$ denote an open neighbourhood of the zero section in
$TM$ as in the introduction. Then we choose a $\mathcal{U}$-admissible
vertical Poisson structure $\theta$ on $TM$ with its corresponding
vertical action of $\mathbb{R}^d$ and (vertical) fundamental vector
fields $X_1, \ldots, X_d$. The corresponding action on the function
spaces is denoted by $\alpha$ as before.
\begin{lemma}
    \label{lemma:CofinalIsometrisch}
    The action $\alpha$ on $C^0(TM)$ is cofinally isometric.
\end{lemma}
\begin{proof}
    Let $K_0 \subseteq M$ be compact. Then for any compact subset $L
    \subseteq TM$ such that $\supp X_i \cap \pi^{-1}(K_0) \subseteq L$
    we have $\norm{f}_L = \norm{\alpha_v(f)}_L$ for all $f \in
    C^0(TM)$ and $v \in \mathbb{R}^d$. Clearly, the compact subsets
    $L$ obeying this condition form a cofinal subset.
\end{proof}

Thus we can define a noncommutative product $\star : C^0(TM)^\infty
\times C^0(TM)^\infty \longrightarrow C^0(TM)^\infty$ by
\begin{equation}
    \label{eq:StarTM}
    f \star g 
    = \iint \alpha_{\Theta u}(f)\; \alpha_v(g) \E^{2\pi \I u \cdot v} 
    \D u\D v.
\end{equation}
Moreover, by Rieffel's construction we also obtain noncommutative
products for $C^0_b(TM)^\infty$ and $C^0_\infty(TM)^\infty$ which we
denote by $\star$, again. Finally, we get noncommutative products for
$C^0(\mathcal{U})^\infty$, $C^0_b(\mathcal{U})^\infty$, and
$C^0_\infty(\mathcal{U})^\infty$ thanks to the support properties of
the action $\alpha$. All these products are continuous for the
deformed (pro-) $C^*$-topologies according to the
Appendices~\ref{subsec:cstar} and \ref{subsec:procstar}.

Next we consider the map $\Phi: TM \supseteq \mathcal{U}
\longrightarrow \mathcal{V}\subseteq M \times M : v_p\mapsto \Phi(v_p)
= (\exp_p(-v_p), \exp_p(v_p))$ that allows to define vector fields
$\tilde{X}_i \in \Gamma^\infty(T(M \times M))$ by
\begin{equation}
    \label{eq:XtildeDef}
    \tilde{X}_i  = \left\{\begin{array}{cl}
            \Phi_* X_i & \textrm{on } \mathcal{V}\\
            0 & \textrm{else.}
        \end{array}\right.
\end{equation}
The corresponding action of $\mathbb{R}^d$ on $C^0(M\times M)$
obtained from these vector fields will be denoted by $\tilde{\alpha}$.
Clearly, $\tilde{\alpha}$ is again cofinally isometric.  Hence this
action gives rise to a deformed product $\tstar : C^0(M \times
M)^\infty \times C^0(M\times M)^\infty \longrightarrow C^0(M\times
M)^\infty$. Analogously to the case of $TM$ we also obtain deformed
products for $C^0_b(M \times M)^\infty$, $C^0_\infty(M \times
M)^\infty$, $C^0(\mathcal{V})^\infty$, $C^0_b(\mathcal{V})^\infty$,
and $C^0_\infty(\mathcal{V})^\infty$ which are continuous for the
corresponding deformed (pro-) $C^*$-topologies.  Now, we relate the
various deformed products:
\begin{proposition}
    \label{proposition:VonTMNachU}
    The restriction map induces a $^*$-homomorphism
    $(C^0(TM)^\infty, \star) \longrightarrow (C^0(\mathcal{U})^\infty,
    \star)$ which is continuous with respect to the deformed
    pro-$C^*$-topologies.  The same statement holds for the case
    $C^0_b$ with respect to the deformed $C^*$-topologies. Moreover,
    the inclusion $(C^0_\infty(\mathcal{U})^\infty, \star)
    \longrightarrow (C^0_\infty(TM)^\infty, \star)$ is a continuous
    $^*$-homomorphism with respect to the deformed $C^*$-topologies.
    Finally, the analogous results hold for $M \times M$ and
    $\mathcal{V}$ instead of $TM$ and $\mathcal{U}$.
\end{proposition}
\begin{proof}
    This is clear from the support property of $\alpha$ and
    Proposition~\ref{proposition:Morphisms}.
\end{proof}

There is no direct relation between the deformed algebras
$(C^0(TM)^\infty, \star)$ and $C^0(M \times M)^\infty, \tstar)$ but we
can relate $\star$ and $\tstar$ restricted to $\mathcal{U}$ and
$\mathcal{V}$, respectively:
\begin{proposition}
    \label{proposition:PhiobenSternEigenschaften}
    The pull-back $\Phi^*: C^0(\mathcal{V}) \longrightarrow
    C^0(\mathcal{U})$ restricts to a $^*$-isomorphism
    \begin{equation}
        \label{eq:PhiObenStern}
        \Phi^*: (C^0(\mathcal{V})^\infty, \tstar) \longrightarrow
        (C^0(\mathcal{U})^\infty, \star)
    \end{equation}
    continuous with respect to the deformed pro-$C^*$-topologies. It
    restricts to a $^*$-isomorphism $\Phi^*:
    (C^0_b(\mathcal{V})^\infty, \tstar) \longrightarrow
    (C^0_b(\mathcal{U})^\infty, \star)$ and $\Phi^*
    :(C^0_\infty(\mathcal{V})^\infty, \tstar) \longrightarrow
    (C^0_\infty(\mathcal{U})^\infty, \star)$, continuous with respect
    to the $C^*$-topolgies. In all cases, the inverse is give by
    $\Phi_*$.
\end{proposition}
\begin{proof}
    First note that $\Phi$ is equivariant with respect to the actions
    $\alpha$ and $\tilde{\alpha}$ whence $\Phi^*$ maps smooth vectors
    to smooth vectors. The Proposition~\ref{proposition:Morphisms} for
    the pro-$C^*$-algebraic situation as well as Rieffel's result
    \cite[Proposition~2.10]{rieffel:1993a} for the $C^*$-algebraic
    case imply the remaining statements.
\end{proof}

In order to obtain a deformed product $\tstarp$ on the functions
$C^0(M)^\infty$ on $M$ that is only noncommutative in a small
neighbourhood of a given point $p \in M$ we have to proceed in two
steps.

First, we consider the embedding $i_p : T_p M \longrightarrow TM$ of
the tangent space at $p \in M$ into the tangent bundle. We want to
show that this map gives rise to a continuous homomorphism of the
algebras $(C^0(TM)^\infty,\star)$ and $(C^0(T_p M)^\infty, \starp)$
via $i_p^*$, where the latter product is obtained as follows: due to
the verticality of the $X_i$ the restrictions of these vector fields
to $T_p M$ define vector fields $X_i^p= X_i|_{T_p M}\in
\Gamma^\infty(T(T_p M))$ and the restrictions of the flows of the
$X_i$ give rise to diffeomorphisms of $T_p M$ which are easily seen to
coincide with the flows of the $X^p_i$. Using these flows we again get
a strongly continuous and cofinally isometric action $\alpha^p$ of
$\mathbb{R}^d$ but now on $C^0(T_p M)$, which can be used to define a
deformed product $\starp$ on $C^0(T_p M)^\infty$ by
\begin{equation}
    f \starp g = \iint
    \alpha^p_{\Theta u}(f)\;\alpha^p_v(g) \E^{2 \pi \I u \cdot v}
    \D u \D v.
\end{equation}
Again, $\starp$ is defined for $C^0_b(T_pM)^\infty$ and
$C^0_\infty(T_pM)^\infty$ as well.
\begin{proposition}
    \label{proposition:iotaStetig}
    The restriction $i_p^*: C^0(TM) \longrightarrow C^0(T_pM)$ induces
    a $^*$-homomorphism
    \begin{equation}
        \label{eq:iHomomorphism}
        i_p^*: (C^0(TM)^\infty, \star)  \longrightarrow
        (C^0(T_pM)^\infty, \starp),
    \end{equation}
    continuous with respect to the deformed pro-$C^*$-topologies.
    Moreover, the analogous statement holds for the $C^*$-algebraic
    cases $C^0_b$ and $C^0_\infty$.
\end{proposition}
\begin{proof}
    The proof is completely analogous to that of
    Proposition~\ref{proposition:PhiobenSternEigenschaften}.
\end{proof}

Now we are prepared to turn to the second step of the construction of
a product for functions on $M$ that is only noncommutative in a small
neighbourhood of $p$. We consider $\mathcal{V}_p = \exp_p
(\mathcal{U}_p)$, where $\mathcal{U}_p = \mathcal{U}\cap T_p M$, and
define vector fields $\tilde{X}_i^p \in \Gamma^\infty (TM)$ on $M$ by
\begin{equation}
    \tilde{X}_i^p = 
    \left\{\begin{array}{cl} 
            (\exp_p)_*X_i^p & \textrm{on }\mathcal{V}_p\\ 
            0 & \textrm{else.}
        \end{array}
    \right.
\end{equation}
Due to the fact that the flows of related vector fields are also
related we get that the composition of the flows of the vector fields
$\tilde{X}_i^p$ coincides with the conjugation of the corresponding
composition of the flows of the $X_i^p$ with $\exp_p$ on
$\mathcal{V}_p$ and extends to an action of $\mathbb{R}^d$ via the
identity outside of $\mathcal{V}_p$. Again we can apply Rieffel's
construction to obtain a product on $C^0(M)^\infty$ as well as on
$C^0_b(M)^\infty$ and $C^0_\infty(M)^\infty$ which will be denoted by
$\tstarp$ and find the following properties.
\begin{proposition}
    Let $p \in M$.
    \begin{enumerate}
    \item $\starp$ restricts from $T_pM$ to $\mathcal{U}_p$ and
        $\tstarp$ restricts from $M$ to $\mathcal{V}_p$, analogously
        to Proposition~\ref{proposition:VonTMNachU}.
    \item The pull-back $\exp_p^*: C^0(\mathcal{V}_p) \longrightarrow
        C^0(\mathcal{U}_p)$ induces a $^*$-isomorphism
        \begin{equation}
            \label{eq:ExppObenStern}
            \exp_p^*: (C^0(\mathcal{V}_p)^\infty, \tstarp) 
            \longrightarrow (C^0(\mathcal{U}_p^\infty), \starp),
        \end{equation}
        continuous with respect to the deformed pro-$C^*$-topologies.
    \item The analogous statements hold for the $C^*$-algebraic cases
        $C^0_b$ and $C^0_\infty$.
    \end{enumerate}
\end{proposition}
\begin{proof}
    The proof is completely analogous to those of
    Proposition~\ref{proposition:VonTMNachU} and
    Proposition~\ref{proposition:PhiobenSternEigenschaften}.
\end{proof}
\begin{remark}
    \label{remark:Frechet}
    Note that all the above $^*$-homomorphisms are also continuous
    with respect to the Fr\'{e}chet topologies of the spaces of smooth
    vectors. This already follows from Rieffel's construction.
\end{remark}
\begin{remark}
    \label{remark:Outlook}
    With these relations we have found the (pro-) $C^*$-algebraic and
    hence non-perturbative counterparts to the constructions from
    \cite{bahns.waldmann:2006a} where formal star products are used
    instead. In order to better understand the physical interpretation
    as indicated in \cite{bahns.waldmann:2006a} one should now
    investigate the state spaces for the deformed algebras as
    explicitly as possible. Moreover, it is interesting to have a
    closer look at the dependence of the deformed products on
    $\theta$. Note that this goes beyond Rieffel's results on
    continuous fields \cite[Chapter~9]{rieffel:1993a} as varying
    $\theta$ implies in particular to vary the action itself.  A good
    understanding of this will be important for interpreting $\theta$
    as a dynamical quantity in more realistic physical models.
    Finally, the results in
    \cite{bahns.doplicher.fredenhagen.piacitelli:2003a,
      doplicher.fredenhagen.roberts:1995a} suggest that one can now
    start developping quantum field theories on locally noncommutative
    spacetimes. One obvious conceptual question is how to interpret
    quantum fields on $TM$.
\end{remark}

%
%

\appendix

%
%

\section{Strict Deformation Quantization for Pro-$C^*$-Algebras}
\label{sec:Appendix}

In \cite{rieffel:1993a}, Rieffel constructs a convergent
noncommutative product for a dense subalgebra of a given Fr{\'e}chet
algebra $\mathcal{A}$. Further he proves, that the resulting
noncommutative algebra is even a pre-$C^*$-algebra, if $\mathcal{A}$
is already a $C^*$-algebra. In this appendix we want to give a short
survey on Rieffel's construction. Moreover, we will show that the
latter result also holds for pro-$C^*$-algebras.

%
%

\subsection{The Deformed Product}
\label{subsec:deformed}

We consider a Fr{\'e}chet algebra $\mathcal{A}$ with a strongly
continuous action $\alpha$ of a vector space $V$ of dimension $d$.
Moreover, $\alpha$ is required to be isometric, i.e.\  there exists a
family $\mathfrak{P}$ of continuous seminorms defining the topology of
$\mathcal{A}$ such that $p(\alpha_v(a)) = p(a)$ for all $p \in
\mathfrak{P}$, $a \in \mathcal{A}$, and $v \in V$.  The space of
smooth vectors for the action $\alpha$ will be denoted by
$\mathcal{A}^\infty$ and is defined by
\begin{equation}
    \mathcal{A}^\infty
    = \left\{ a \in \mathcal{A} 
        \Bigg| \left. 
            \left(\frac{\D}{\D t_1}\right)^{\gamma_1} \ldots
            \left(\frac{\D}{\D t_d}\right)^{\gamma_d} \alpha_{\exp(t_1 e_1)\cdots
              \exp(t_d e_d)}(a)\right|_{t_1=\ldots=t_d=0} 
        \; \textrm{exists for all} \; \gamma \in \mathbb{N}_0^d \right\},
\end{equation}
where $e_1, \ldots, e_d$ form a basis of $V$.  $\mathcal{A}^\infty$
carries a Fr{\'e}chet topology in an obvious way.  Moreover, it is
dense in $\mathcal{A}$ with respect to the original Fr{\'e}chet
topology \cite[Theorem A.1]{schweitzer:1993a}.  Further, let $\tau$ be
the action of $V$ on the space of uniformly continuous bounded
mappings from $V$ to $\mathcal{A}$, $C^0_u(V,\mathcal{A})$, by
translation. With a convenient partition of unity $\{\varphi_w\}$ for
$V \times V$, Rieffel shows that the oscillating integrals
\begin{equation}
    \iint F(u,v) \E^{2 \pi \I u\cdot v} \D u \D v 
    = \sum\limits_w \iint
    F(u,v) \varphi_w (u,v) \E^{2 \pi \I u\cdot v}\D u \D v,
\end{equation}
are well-defined \cite[Proposition 1.6]{rieffel:1993a}.  Here $F$ is
in the space of smooth vectors of the canonical action of $V \times V$
on $C^0_u(V\times V, \mathcal{A})$ by translations.  Then, if $\Theta$
is a skew-symmetric operator on $V$,
\begin{equation}
    a \star b = \iint \alpha_{\Theta u}(a)\alpha_v(b) \E^{2 \pi \I u\cdot v} 
    \D u \D v
\end{equation}
is a noncommutative associative product for $\mathcal{A}^\infty$,
called the {\em deformed product} (determined by $\alpha$ and
$\Theta$), see \cite[Definition 2.1, Theorem 2.14]{rieffel:1993a}.
Rieffel shows various properties of this deformed product among of
which we need the following \cite[Corollary 2.13]{rieffel:1993a}:
\begin{proposition}
    \label{proposition:fixpoint}
    Let $a \in \mathcal{A}^\infty$ be a fixed point for the action
    $\alpha$. Then for any $b \in \mathcal{A}^\infty$, we have $a \star b
    = ab$ and $b \star a = ba$. 
\end{proposition}
\begin{remark}
    \label{remark:JenseitsFrechet}
    Note that the construction so far can be carried through also for
    a Hausdorff complete locally convex topological algebra not
    necessarily Fr{\'e}chet.
\end{remark}

%
%

\subsection{$C^*$-Algebras}
\label{subsec:cstar}

Let $\mathcal{A}$ be a $C^*$-algebra. In order to define a $C^*$-norm
for the deformed algebra $(\mathcal{A}^\infty, \star )$, one considers
the space of all functions in $C^0_u(V, \mathcal{A})^\infty$ such that
the product of their derivatives with any polynomial on $V$ are
bounded which will be denoted by ${\mathcal{S^A}(V)}$. On
${\mathcal{S^A}(V)}$, an $\mathcal{A}$-valued inner product is defined
by
\begin{equation}
    \langle {f}, {g} \rangle = \int f(v)^* g(v) \D v.
\end{equation}
As in the case of Hilbert spaces, a corresponding norm for
${\mathcal{S^A}(V)}$ is defined by:
\begin{equation}
    \norm{f} = \norm{\langle {f}, {f} \rangle }^{\frac{1}{2}}.
\end{equation}
Since $C^0_u(V, \mathcal{A})^\infty$ carries the action $\tau$ of $V$
and is itself a Fr{\'e}chet algebra, we obtain a deformed product
$\star$ for $C^0_u(V, \mathcal{A})^\infty$.  Furthermore, let
$\mathsf{L}$ be the action of $C^0_u(V, \mathcal{A})^\infty$ on
${\mathcal{S^A}(V)}$ given by this deformed product
\begin{equation}
    \mathsf{L}_Fg = F \star g.
\end{equation}
Rieffel shows, that the operator $\mathsf{L}_f: {\mathcal{S^A}(V)}
\longrightarrow {\mathcal{S^A}(V)}$ is bounded and adjointable for $f
\in {\mathcal{S^A}(V)}$, see \cite[Corollary 4.4]{rieffel:1993a}.
Further, he proves that this result also holds for the case $F \in
C^0_u(V, \mathcal{A})^\infty$, see \cite[Theorem~4.6]{rieffel:1993a}.

Let $\alpha$ be again an isometric strongly continuous action of $V$
on the $C^*$-algebra $\mathcal{A}$. For any $a \in \mathcal{A}$, we
define a function $\phi(a)$ by
\begin{equation}
    \phi(a)(v) = \alpha_v(a).
\end{equation}
Since $\alpha$ is isometric it turns out that $\phi(a) \in C_u^0(V,
\mathcal{A})$.  The map $a \mapsto \phi(a)$ constitutes a
$^*$-homomorphism from $\mathcal{A}$ into $C^0_u(V, \mathcal{A})$,
which is equivariant with respect to the actions $\alpha$ on
$\mathcal{A}$ and $\tau$ on $C^0_u(V, \mathcal{A})$. Thus it carries
smooth vectors to smooth vectors, i.e.\ $\mathcal{A}^\infty$ into
$C^0_u(V, \mathcal{A})^\infty$, and is a homomorphism for their
deformed products \cite[Proposition 2.10]{rieffel:1993a}. Therefore,
each $a \in \mathcal{A}^\infty$ determines a bounded operator on
${\mathcal{S^A}(V)}$ which will be denoted by $\mathsf{L}_{\phi(a)}$.
Thus we can define a norm on $\mathcal{A}^\infty$ by
\begin{equation}
    \label{eq:NormTheta}
    \norm{a}_\Theta = \norm{\mathsf{L}_{\phi(a)}}.
\end{equation} 
As the adjointable operators on a Hilbert module over a $C^*$-algebra
form again a $C^*$-algebra \cite[p.~8]{lance:1995a}, it is obvious
that the defined norm $\norm{\,.\,}_\Theta$ satisfies the
$C^*$-conditions. Thus the completion of $(\mathcal{A}^\infty, \star,
\norm{\,.\,}_\Theta)$ is indeed a $C^*$-algebra. We shall refer to
this $C^*$-topology as the deformed topology.

%
%

\subsection{Pro-$C^*$-Algebras}
\label{subsec:procstar}

We shall now extend Rieffel's construction to the case of
pro-$C^*$-algebras.  Let $\mathcal{A} = \varprojlim
\mathcal{A}_\lambda$ be a pro-$C^*$-algebra, i.e.\ the inverse limit
of an inverse system of $C^*$-algebras $\mathcal{A}_\lambda$ in the
category of topological $^*$-algebras, see e.g.~\cite{phillips:1988a,
  phillips:1988b}. An inverse system of $C^*$-algebras consists of a
directed set $\Lambda$, a $C^*$-algebra $\mathcal{A}_\lambda$ for each
$\lambda \in \Lambda$, and $^*$-homomorphisms $\pi_{\lambda, \rho}:
\mathcal{A}_\lambda \longrightarrow \mathcal{A}_\rho$ for $\lambda \ge
\rho$, satisfying the following conditions:
\begin{equation}
    \label{eq:ProjectiveLimit}
    \pi_{\lambda, \lambda} 
    = \id_{\mathcal{A}_\lambda} 
    \quad
    \textrm{and}
    \quad
    \pi_{\rho, \mu} \circ \pi_{\lambda, \rho} 
    = \pi_{\lambda, \mu} 
    \quad \textrm{for} \; \lambda \ge \rho \ge \mu.
\end{equation}
The inverse limit of the system ($\mathcal{A}_\lambda, \pi_{\lambda,
  \rho}$) in the category of topological $^*$-algebras is a
topological $^*$-algebra $\mathcal{A}$ together with
$^*$-homomorphisms $\kappa_\lambda: \mathcal{A} \longrightarrow
\mathcal{A}_\lambda$, such that
\begin{equation*}
    \pi_{\lambda, \mu} \circ \kappa_\lambda = \kappa_\mu
\end{equation*}
and satisfying the usual universal property as in
\cite{phillips:1988a}.  An element $a \in \mathcal{A}$ can be
identified with a coherent sequence $(a_\lambda) \in
\prod\limits_{\lambda \in \Lambda} \mathcal{A}_\lambda$ satisfying
$\pi_{\lambda, \mu} a_\lambda = a_\mu$. The topology of $\mathcal{A}$
is determined by the set of all continuous $C^*$-seminorms on
$\mathcal{A}$ denoted by $S(\mathcal{A})$, see
\cite[Proposition~1.1.1]{phillips:1988a}. $S(\mathcal{A})$ is
obviously a directed set by $q' \ge q$ iff $q'(a) \ge q(a)$ for all $a
\in \mathcal{A}$.  Defining $\mathcal{A}_q = \mathcal{A} / \ker q$ for
$q \in S(\mathcal{A})$ one has by
\cite[Proposition~1.1.1]{phillips:1988a}
\begin{equation}
    \mathcal{A} \cong \varprojlim \mathcal{A}_q.
\end{equation}

An important example for a pro-$C^*$-algebra we use in
Section~\ref{sec:RieffelCompactSupport} is the algebra of continuous
functions $C^0(N)$ over a manifold $N$. It is the inverse limit of the
inverse system of the $C^*$-algebras $(C^0(L), \norm{\,.\,}_L)$ where
\begin{equation*}
    \norm{f}_L = \sup_{x \in L} \left|{f(x)}\right|.
\end{equation*}
The compact sets are ordered by $L' \ge L$ iff $L' \supseteq L$.  The
mappings $\pi_{L', L}: C^0(L') \longrightarrow C^0(L)$ for $L \subseteq
L'$ are given by
\begin{equation*}
    \pi_{L',L} f = i^*_L f,
\end{equation*}
where $i^*_L$ is pull-back to $L$. The mappings $\kappa_L: C^0(N)
\longrightarrow C^0(L)$ are given by
\begin{equation*}
    \kappa_{L} f = i^*_L f
\end{equation*}

For a strongly continuous action $\alpha$ of $V$ on a
pro-$C^*$-algebra $\mathcal{A}$ we introduce the following definition:
\begin{definition}
    \label{definition:CofinaleWirkung}
    The action $\alpha$ is called cofinally isometric if there exists
    a cofinal subset $\Lambda \subseteq S(\mathcal{A})$ such that for
    all $q \in \Lambda$, all $a \in \mathcal{A}$, and all $v \in V$ we
    have
    \begin{equation}
        \label{eq:CofinalIsometrisch}
        q(\alpha_v(a)) = q(a).
    \end{equation}
\end{definition}

Given such a cofinally isometric and strongly continuous action
$\alpha$ on $\mathcal{A}$, we obtain a deformed product $\star$ on the
smooth vectors $\mathcal{A}^\infty$ by the general results of
Section~\ref{subsec:deformed}. Now we want to define corresponding
$C^*$-seminorms for $(\mathcal{A}^\infty, \star)$ such that the
completion with respect to these seminorms gives again a
pro-$C^*$-algebra.

We can proceed analogously to the case of $C^*$-algebras up to the
definition of the deformed norm as in \eqref{eq:NormTheta}.  Here we
have to be more specific.

Recall that a Hilbert module $\mathcal{E}$ over a pro-$C^*$-algebra
$\mathcal{A}$ is defined analogously to the $C^*$-algebraic case, see
\cite{phillips:1988b}, where completeness is now understood with
respect to the seminorms $\norm{\xi}_q = q(\SP{\xi,
  \xi})^{\frac{1}{2}}$ where $\xi \in \mathcal{E}$ and $q \in
S(\mathcal{A})$. Then $\mathcal{E}_q$ is defined to be the quotient
$\mathcal{E} \big/ \ker \norm{\,.\,}_q$ and turns out to be a Hilbert
module over the $C^*$-algebra $\mathcal{A}_q$. Thus the continuous
adjointable operators $\mathfrak{B}(\mathcal{E}_q)$ on $\mathcal{E}_q$
form a $C^*$-algebra with respect to the usual operator norm.  Given a
continuous adjointable operator $T \in \mathfrak{B}(\mathcal{E})$, one
defines $T_q \in \mathfrak{B}(\mathcal{E}_q)$ by $T_q [\xi]_q =
[T\xi]_q$, where $[\xi]_q \in \mathcal{E}_q$ denotes the class of
$\xi$. Then $\mathfrak{B}(\mathcal{E}) \ni T \mapsto T_q \in
\mathfrak{B}(\mathcal{E}_q)$ is clearly a $^*$-homomorphism whence
\begin{equation}
    \label{eq:NormTq}
    \norm{T}_q = \norm{T_q}
\end{equation}
defines a $C^*$-seminorm for $\mathfrak{B}(\mathcal{E})$ for each $q
\in S(\mathcal{A})$. The pro-$C^*$-topology induced by these seminorms
coincides with the one in \cite{phillips:1988b}.

Using this pro-$C^*$-topology for
$\mathfrak{B}(\mathcal{S}^{\mathcal{A}}(V))$ we can define the
pro-$C^*$-seminorms
\begin{equation}
    \label{eq:ProCstarSeminorm}
    \norm{a}_{\Theta, q} = \norm{\mathsf{L}_{\phi(a)}}_q
\end{equation}
for $a \in (\mathcal{A}^\infty, \star)$. The completion of
$(\mathcal{A}^\infty, \star)$ with respect to this deformed
pro-$C^*$-topology completes the construction.
\begin{proposition}
    \label{proposition:ProCStarDeformation}
    For a pro-$C^*$-algebra $\mathcal{A}$ endowed with a strongly
    continuous and cofinally isometric action $\alpha$ of $V$ the
    deformed algebra $(\mathcal{A}^\infty, \star)$ of smooth vectors
    carries a system of $C^*$-seminorms defined by
    \eqref{eq:ProCstarSeminorm}. The completion with respect to these
    seminorms yields a pro-$C^*$-algebra deforming $\mathcal{A}$.
\end{proposition}

We conclude with a last remark on the functoriality of the
construction: Let $\mathcal{A}$ and $\mathcal{B}$ be
pro-$C^*$-algebras and $\Psi: \mathcal{A} \longrightarrow \mathcal{B}$
be a continuous $^*$-homomorphism. Suppose that $\mathcal{A}$ and
$\mathcal{B}$ are equipped with strongly continuous and cofinally
isometric actions of $V$ such that $\Psi$ is equivariant.
\begin{proposition}
    \label{proposition:Morphisms}
    With the assumptions from above $\Psi: (\mathcal{A}^\infty, \star)
    \longrightarrow (\mathcal{B}^\infty, \star)$ is a continuous
    $^*$-homomorphism with respect to the deformed
    pro-$C^*$-topologies.
\end{proposition}
\begin{proof}
    First, one shows that $\kappa_p: \mathcal{A} \longrightarrow
    \mathcal{A}_p$ induces a $^*$-homomorphism $\kappa_p:
    \mathcal{A}^\infty \longrightarrow \mathcal{A}^\infty_p$ which
    turns out to be continuous with respect to the deformed (pro-)
    $C^*$-topologies. In fact, for $a \in \mathcal{A}^\infty$ we have
    $\norm{\kappa_p(a)}_{\Theta} = \norm{a}_{\Theta, p}$ by a
    straightforward computation. This implies that the projective
    limit of the deformations $(\mathcal{A}_p^\infty, \star_p,
    \norm{\,.\,}_\Theta)$ is isomorphic to the deformation
    $(\mathcal{A}^\infty, \star, \{\norm{\,.\,}_{\Theta,p}\}_{p \in
      S(\mathcal{A})})$ of the projective limit $\mathcal{A}$.
    Clearly, the same is true for $\mathcal{B}$.  Second, since $\Psi$
    is continuous with respect to the undeformed pro-$C^*$-topologies,
    we find for each $q \in S(\mathcal{B})$ a $p \in S(\mathcal{A})$
    such that $q(\Psi(a)) = p(a)$ for all $a \in \mathcal{A}$. This
    implies that for these $q$, $p$ we obtain a $^*$-homomorphism
    $\Psi_{p,q}: \mathcal{A}_p \longrightarrow \mathcal{B}_q$ with
    respect to the undeformed $C^*$-topologies. Since $\Psi$ is
    equivariant, it gives a continuous $^*$-homomorphism $\Psi_{p,q}:
    \mathcal{A}^\infty_p \longrightarrow \mathcal{B}^\infty_q$ with
    respect to the deformed $C^*$-topologies according to
    \cite[Theorem~5.7]{rieffel:1993a}.  Thus we have continuous
    $^*$-homomorphisms $\Psi_q = \Psi_{p,q} \circ \kappa_p = \kappa_q
    \circ \Psi: \mathcal{A}^\infty \longrightarrow
    \mathcal{B}^\infty_q$ which, by the universal property of
    projective limits, after completion combine to a continuous
    $^*$-homomorphism between the deformed pro-$C^*$-algebras.
    Clearly, on $\mathcal{A}^\infty$ it coincides with $\Psi$.
\end{proof}

%
%

\begin{footnotesize}
    \renewcommand{\arraystretch}{0.5} 

\end{footnotesize}
\end{document}